\documentclass[12pt,reqno,letter]{amsart}

\usepackage[all]{xy}
\usepackage{texdraw}

\def\theoremnumberstyle{plain}
\def\proofname{Proof}
\usepackage{tom}

\newcommand{\tom}[1]{}

\newcommand{\pp}{{\mathbb P}}

\SelectTips{cm}{12}

\begin{document}

   \parindent0cm

   \title{A Note on Equimultiple Deformations}
   \author{Thomas Markwig}
   \address{Universit\"at Kaiserslautern\\
     Fachbereich Mathematik\\
     Erwin-Schr\"odinger-Stra\ss e\\
     D -- 67663 Kaiserslautern
     }
   \email{keilen@mathematik.uni-kl.de}
   \urladdr{http://www.mathematik.uni-kl.de/\textasciitilde keilen}
   \thanks{The author was supported by the EAGER node of Torino.}


   \date{October, 2006.}


   \begin{abstract}
     While the tangent space to an equisingular family of curves can be
     discribed by the sections of a twisted ideal sheaf, this is no
     longer true if we only prescribe the multiplicity which a singular
     point should have. However, it is still possible to compute the
     dimension of the tangent space with the aid of the
     equimulitplicity ideal.
     In this note we consider families $\kl_m=\{(C,p)\in|L|\times
     S\;|\;\mult_p(C)=m\}$ for some linear system $|L|$ on a smooth
     projective surface $S$ and a fixed positive integer $m$, and we
     compute the dimension of the tangent space to $\kl_m$ at a point
     $(C,p)$ depending on whether $p$ is a unitangential singular
     point of $C$ or not. We deduce that the expected dimension of
     $\kl_m$ at $(C,p)$ in any case is just
     $\dim|L|-\frac{m\cdot(m+1)}{2}+2$. The result is used in the
     study of triple-point defective surfaces in \cite{CM07} and
     \cite{CM07a}.  
   \end{abstract}

   \maketitle

   The paper is based on considerations about the Hilbert scheme of
   curves in a projective surface (see e.g.\ \cite{Mum66},
   Lecture 22) and about local equimultiple deformations of plane
   curves (see \cite{Wah74}).

   \begin{definition}
     Let $T$ be a complex space. An \emph{embedded family of curves in
       $S$ with section} over $T$ is a commutative diagram of morphisms
     \begin{displaymath}
       \xymatrix{
         \kc\ar@{^{(}->}[r]\ar[d]_\varphi&T\times S\ar[dl]\\
         T\ar@(l,l)[u]^\sigma
         }
     \end{displaymath}
     where $\codim_{T\times S}(\kc)=1$, $\varphi$ is flat and proper,
     and $\sigma$ is a section, i.e.\
     $\varphi\circ\sigma=\id_T$. Thus we have a morphism
     $\ko_T\rightarrow\varphi_*\ko_\kc=\varphi_*\big(\ko_{T\times S}/\kj_\kc\big)$
     such that $\varphi_*\ko_\kc$ is a flat $\ko_T$-module.

     The family is said to be \emph{equimultiple} of multiplicity $m$
     \emph{along the section} $\sigma$ if the ideal sheaf $\kj_\kc$  of $\kc$ in
     $\ko_{T\times S}$ satisfies
     \begin{displaymath}
       \kj_\kc\subseteq\kj_{\sigma(T)}^m\;\;\;\mbox{ and }
       \;\;\;\kj_\kc\not\subseteq\kj_{\sigma(T)}^{m+1},
     \end{displaymath}
     where $\kj_{\sigma(T)}$ is the ideal sheaf  of $\sigma(T)$ in
     $\ko_{T\times S}$.
   \end{definition}

   \begin{remark}
       Note that the above notion commutes with base change, i.e.\ if we
       have an equimultiple embedded family of curves in $S$ over $T$ as
       above and if $\alpha:T'\rightarrow T$ is a morphism, then the fibre
       product diagram
       \begin{displaymath}
         \xymatrix{
           T'\times S\ar[rrrr]\ar[ddr]&&&& T\times S\ar[ddl]\\
           &\kc'\ar@{_{(}->}[ul]\ar[d]^{\varphi'}\ar[rr]&&\kc\ar@{^{(}->}[ur]\ar[d]_\varphi\\
           &T'\ar@(r,r)[u]^{\sigma'}\ar[rr]^\alpha&&T\ar@(l,l)[u]_\sigma
         }
       \end{displaymath}
       gives rise to an embedded equimultiple family of curves over $T'$
       of the same multiplicity, since locally it is defined via the
       tensor product.
   \end{remark}

   \begin{example}\label{ex:infinitesimal}
     Let us denote by $T_\varepsilon=\Spec(\C[\varepsilon])$ with
     $\varepsilon^2=0$. Then a family of curves in $S$ over
     $T_\varepsilon$ is just a Cartier divisor of $T_\varepsilon\times
     S$, that is, it is given on a suitable open covering
     $S=\bigcup_{\lambda\in\Lambda}U_\lambda$ by equations
     \begin{displaymath}
       f_\lambda+\varepsilon\cdot g_\lambda\in
       \C[\varepsilon]\otimes_\C
       \Gamma(U_\lambda,\ko_S)=\Gamma(U_\lambda,\ko_{T\times S}),
     \end{displaymath}
     which glue together to give a global section
     $\big\{\frac{g_\lambda}{f_\lambda}\big\}_{\lambda\in\Lambda}$
     in $H^0\big(C,\ko_C(C)\big)$, where $C$ is the curve defined  locally by the
     $f_\lambda$ (see e.g.\ \cite{Mum66}, Lecture 22).
     \tom{
       And vice versa, every section $s\in
       H^0\big(S,\ko_S(C)\big)$ for a suitable covering
       $S=\bigcup_{\lambda\in\Lambda}U_\lambda$  is given as
       \begin{displaymath}
         s_{|U_\lambda}=\frac{g_\lambda}{f_\lambda}
       \end{displaymath}
       where $f_\lambda,g_\lambda\in\Gamma(U_\lambda,\ko_S)$ are regular
       functions on $U_\lambda$ such that $f_\lambda$ is a local equation
       for $C\cap U_\lambda$ and
       $\forall\;\lambda,\lambda'\in\Lambda\;:\;\frac{g_\lambda}{f_\lambda}\equiv
       \frac{g_{\lambda'}}{f_{\lambda'}}$ on $U_\lambda\cap
       U_{\lambda'}$, i.e.\
       $\frac{f_\lambda}{f_{\lambda'}}
       \in\Gamma(U_\lambda\cap
       U_{\lambda'},\ko_S)^*$   and
       $g_\lambda=\frac{f_\lambda}{f_{\lambda'}}\cdot g_{\lambda'}$ in
       $\Gamma(U_\lambda\cap U_{\lambda'},\ko_S)$. Thus the
       $f_\lambda+\varepsilon\cdot g_\lambda$ are the equations of a
       first-order infinitesimal deformation of $C$ in $|L|$.
     }

     A section of the family through $p$ is locally in $p$ given as
     $(x,y)\mapsto(x_a,y_b)=(x+\varepsilon\cdot a,y+\varepsilon\cdot
     b)$ for some $a,b\in\C\{x,y\}=\ko_{S,p}$.
   \end{example}

   \begin{example}\label{ex:incidencevariety}
     Let  $H$ be a connected component of the Hilbert scheme
     $\Hilb_S$ of curves in $S$, then $H$ comes with a
     universal family
     \begin{equation}\label{eq:def:1}
       \pi:\kh\longrightarrow H:(C,p)\mapsto C.
     \end{equation}
     Let us now fix a positive integer $m$ and set
     \begin{displaymath}
       {\kh_m}=\{(C,p)\in H\times S\;|\;C\in H,\mult_p(C)=m\}.
     \end{displaymath}
     Then ${\kh_m}$ is a locally closed subvariety of $H\times S$, and
     \eqref{eq:def:1} induces via base change a flat and proper family ${\kf_m}=\{(C_p,q)\in
     {\kh_m}\times S\;|\;
     C_p=(C,p)\in {\kh_m}, q\in C\}$ which has a
     distinguished section $\sigma$
     \begin{equation}\label{eq:def:1a}
       \xymatrix{
         {\kf_m}\;\;\ar@{^{(}->}[r]\ar[d]
         &\;\;{\kh_m}\times S\ar[dl]\\
         {\kh_m}\ar@(l,l)[u]_\sigma
       }
     \end{equation}
     sending $C_p=(C,p)$ to $(C_p,p)\in{\kf_m}$. Moreover, this family is
     equimultiple along $\sigma$ of multiplicity $m$ by construction.
%
   \end{example}

   \begin{example}\label{ex:incidencevariety-ls}
     Similarly, if $|L|$ is a linear system on $S$,  then it induces a
     universal family
     \begin{equation}\label{eq:def:1ls}
       \pi:\kl=\{(C,p)\in|L|\times S\;|\;p\in C\}\longrightarrow
       |L|:(C,p)\mapsto C.
     \end{equation}
     If we now fix a positive integer $m$ and set
     \begin{displaymath}
       {\kl_m}=\{(C,p)\in|L|\times S\;|\;C\in |L|,\mult_p(C)=m\}.
     \end{displaymath}
     Then ${\kl_m}$ is a locally closed subvariety of $|L|\times S$, and
     \eqref{eq:def:1ls} induces via base change a flat and proper family $\kg_m=\{(C_p,q)\in
     {\kl_m}\times S\;|\;
     C_p=(C,p)\in {\kl_m}, q\in C\}$ which has a
     distinguished section $\sigma$
     \begin{equation}\label{eq:def:1als}
       \xymatrix{
         \kg_m\;\;\ar@{^{(}->}[r]\ar[d]
         &\;\;{\kl_m}\times S\ar[dl]\\
         {\kl_m}\ar@(l,l)[u]_\sigma
       }
     \end{equation}
     sending $C_p=(C,p)$ to $(C_p,p)\in\kg_m$. Moreover, this family is
     equimultiple along $\sigma$ of multiplicity $m$ by construction.

     We may interpret ${\kl_m}$ as the family of curves in $|L|$ with $m$-fold
     points together with a section which distinguishes the $m$-fold point.
     This is important if the $m$-fold point is not isolated or if it
     splits in a neighbourhood into several simpler $m$-fold points.

     Of course, since \eqref{eq:def:1ls} can be viewed as a subfamily
     of \eqref{eq:def:1} we may view \eqref{eq:def:1als} in the same
     way as a subfamily of \eqref{eq:def:1a}.
   \end{example}

   \begin{definition}
     Let $t_0\in T$ be a pointed complex space, $C\subset S$ a curve,
     and $p\in C$ a point of multiplicity $m$. Then an \emph{embedded
       (equimultiple) deformation of $C$ in $S$ over $t_0\in T$ with section $\sigma$ through}
     $p$ is a commutative diagram of morphisms
     \begin{displaymath}
       \xymatrix{
         &S\ar@{^{(}->}[rr]\ar[ddl]|\hole&& T\times S\ar[ddl]\\
         p\in C\;\;\;\;\;\ar@{^{(}->}[ur]\ar[d]\ar@{^{(}->}[rr]&&\kc\ar@{^{(}->}[ur]\ar[d]_\varphi\\
         t_0\ar@(l,l)[u]_\sigma\ar@{|->}[rr]&&T\ar@(ul,dl)[u]^\sigma
       }
     \end{displaymath}
     where the right hand part of the diagram is an embedded
     (equimultiple) family of curves in $S$ over $T$ with section
     $\sigma$. Sometimes we will simply write $(\varphi,\sigma)$ to
     denote a deformation as above.

     Given two deformations, say $(\varphi,\sigma)$ and
     $(\varphi',\sigma')$, of $C$ over $t_0\in T$ as above,
     a morphism of these deformations is a morphism $\psi:\kc'\rightarrow
     \kc$ which makes the obvious diagram commute:
     \begin{displaymath}
            \xymatrix{
       T\times S\ar@{=}[rrrr]\ar[ddddr]&&&& T\times S\ar[ddddl]\\
       &\kc'\ar@{_{(}->}[ul]\ar[ddd]_{\varphi'}\ar[rr]^{\psi}&&\kc\ar@{^{(}->}[ur]\ar[ddd]^\varphi\\
       &&C\ar@{_{(}->}[ul]\ar@{^{(}->}[ur]\ar[d]\\
       &&t_0\ar@{|->}[dr]\ar@{|->}[dl]\\
       &T\ar@(ur,dr)[uuu]^{\sigma'}\ar@{=}[rr]&&T\ar@(ul,dl)[uuu]_\sigma
       }
     \end{displaymath}
     This gives rise to the deformation functor
     \begin{displaymath}
       \underline{\Def}^{sec,em}_{p\in C/S}:(\mbox{pointed complex
         spaces})
       \rightarrow(\mbox{sets})
     \end{displaymath}
     of embedded equimultiple deformations of $C$ with section through
     $p$ from
     the category of pointed complex spaces into the category of sets, where
     for a pointed complex space $t_0\in T$
     \begin{align*}
       \underline{\Def}^{sec,em}_{p\in C/S}(t_0\in T)=\{&
       \mbox{isomorphism classes of embedded equimultiple}
       \\
       &\mbox{deformations $(\varphi,\sigma)$
         of $C$ in $S$ over $t_0\in T$ }
       \\
       &\mbox{with section through $p$}\}.
     \end{align*}
     Moreover, forgetting the section we have a natural transformation
     \begin{equation}\label{eq:def:2}
       \underline{\Def}^{sec,em}_{p\in C/\Sigma}
       \longrightarrow
       \underline{\Def}_{C/\Sigma},
     \end{equation}
     where the latter is the deformation functor
     \begin{displaymath}
       \underline{\Def}_{C/\Sigma}:(\mbox{pointed complex spaces})\rightarrow(\mbox{sets})
     \end{displaymath}
     of embedded deformations of $C$ in $S$ given by
     \begin{align*}
       \underline{\Def}_{C/S}(t_0\in T)=\{&
       \mbox{isomorphism classes of embedded deformations}
       \\
       &\mbox{of $C$ in $S$ over $t_0\in T$ }\}.
     \end{align*}
   \end{definition}

   \begin{example}
     According to Example \ref{ex:infinitesimal} a deformation of $C$
     in $S$ over $T_\varepsilon$ along a section
     through $p$ is given by
     \begin{itemize}
     \item local equations
       $f+\varepsilon\cdot g$
       such that $f$ is a local equation for $C$ and the $\frac{g}{f}$
       glue to a global section of
       $\ko_C(C)$,
     \item  together with a section which in local coordinates in
       $p$ is given as $\sigma:(x,y)\mapsto
       (x_a,y_b)=(x+\varepsilon\cdot a,y+\varepsilon\cdot b)$ for some
       $a,b\in\C\{x,y\}$.
     \end{itemize}
     If we forget the section it is well known (see e.g.\ \cite{Mum66},
     Lecture 22) that two such
     deformations are isomorphic if and only if they induce the same
     global section of $\ko_C(C)$ and this one-to-one correspondence
     is functorial so that we have an isomorphism of vector spaces
     \begin{displaymath}
       \underline{\Def}_{C/S}(T_\varepsilon)\stackrel{\cong}{\longrightarrow}
       H^0\big(C,\ko_C(C)\big).
     \end{displaymath}

     Considering the natural transformation from \eqref{eq:def:2} we
     may now ask what the image of
     $\underline{\Def}^{sec,em}_{p\in C/S}(T_\varepsilon)$ in
     $H^0\big(C,\ko_C(C)\big)$ is. These are, of course, the sections
     which allow a section $\sigma$ through $p$ along which the deformation is
     equimultiple, and according to Lemma \ref{lem:equimultiple} we
     thus have an epimorphism
     \begin{displaymath}
       \xymatrix{
         \underline{\Def}^{sec,em}_{p\in C/S}(T_\varepsilon)
         \ar@{->>}[r]
         & H^0\big(C,\kj_{Z/C}(C)\big),
       }
     \end{displaymath}
     where $\kj_{Z/C}$ is the restriction to $C$ of the ideal sheaf
     $\kj_Z$ on $S$ given by
     \begin{equation}\label{eq:Z}
       \kj_{Z,q}=\left\{
         \begin{array}{ll}
           \ko_{S,q},& \mbox{ if }q\not=p,\\
           \big\langle\frac{\partial f}{\partial x}, \frac{\partial
               f}{\partial y}\big\rangle+
           \langle x,y\rangle^m, &
           \mbox{ if }q=p,
         \end{array}
       \right.
     \end{equation}
     here $f$ is a local equation for $C$ in local coordinates $x$ and
     $y$ in $p$.

     It remains the question what the dimension of the kernel of this map is, that is,
     how many different sections such an isomorphism class of
     embedded deformations of $C$ in $S$ over $T_\varepsilon$ through
     $p$ can admit.

     J.\ Wahl showed in \cite{Wah74}, Proposition 1.9, that locally the
     equimultiple deformation admits a unique section if and only if $C$
     in $p$ is not unitangential. If $C$ is unitangential we may
     assume that locally in $p$ it is given by $f=y^m+h.o.t.$. If we have an
     embedded deformation of $C$ in $S$ which along some section is
     equimultiple of multiplicity $m$, then locally it looks like
     \begin{displaymath}
       f+\varepsilon\cdot \left(a\cdot {\textstyle\frac{\partial f}{\partial x}}
         +b\cdot {\textstyle\frac{\partial f}{\partial y}}+h\right)
     \end{displaymath}
     with $h\in \langle x,y\rangle^m$. However, since $\frac{\partial
       f}{\partial x}\in \langle x,y\rangle^m$ the deformation is
     equimultiple along the sections $(x,y)\mapsto (x+
     \varepsilon\cdot (c+a),y+\varepsilon\cdot b)$ for all $c\in \C$. Thus
     in this case the kernel turns out to be one-dimensional, i.e.\
     there is a one-dimensional vector space $\kk$ such that the
     following sequence is exact:
     \begin{equation}
       \label{eq:def:3}
       0\rightarrow\kk\rightarrow
       \underline{\Def}^{sec,em}_{p\in C/S}(T_\varepsilon)
       \rightarrow
       H^0\big(C,\kj_{Z/C}(C)\big)
       \rightarrow 0.
     \end{equation}
   \end{example}

   \begin{lemma}\label{lem:equimultiple}
     Let $f+\varepsilon\cdot g$ be a first-order infinitesimal deformation  of
     $f\in\C\{x,y\}$, $m=\ord(f)$, $a,b\in\C\{x,y\}$, and
     $x_a=x+\varepsilon\cdot a$, $y_b=y+\varepsilon\cdot b$.

     Then $f+\varepsilon\cdot g$ is equimultiple along the section
     $(x,y)\mapsto(x_a,y_b)$ if and only if
     \begin{displaymath}
       g-a\cdot \frac{\partial f}{\partial x}-b\cdot
       \frac{\partial f}{\partial y}
       \in
       \langle x,y\rangle^m.
     \end{displaymath}
     In particular, $f+\varepsilon\cdot g$ is equimultiple along some
     section if and only if
     \begin{displaymath}
       g\in
       \left\langle\frac{\partial f}{\partial x}, \frac{\partial
           f}{\partial y}\right\rangle+
       \langle x,y\rangle^m.
     \end{displaymath}
   \end{lemma}
   \begin{proof}
     If $a,b\in \C\{x,y\}$ and $h\in\langle
     x,y\rangle^m$ then by Taylor expansion and since $\varepsilon^2=0$ we have
     \begin{displaymath}
       f+\varepsilon\cdot \left(a\cdot {\textstyle\frac{\partial f}{\partial x}}
         +b\cdot {\textstyle\frac{\partial f}{\partial y}}+h\right)
       =f(x_a,y_b)+\varepsilon\cdot h(x_a,y_b),
     \end{displaymath}
     where $f(x_a,y_b), h(x_a,y_b)\in \langle
     x_a,y_b\rangle^m$, i.e.\ the infinitesimal deformation $f+\varepsilon\cdot
     \big(a\cdot {\textstyle\frac{\partial f}{\partial x}}
       +b\cdot {\textstyle\frac{\partial f}{\partial y}}+h\big)$
     is equimultiple along $(x,y)\mapsto(x_a,y_b)$.

     Conversely, if $f+\varepsilon \cdot g$ is equimultiple along
     $(x,y)\mapsto(x_a,y_b)$ then
     \begin{displaymath}
       f(x,y)+\varepsilon \cdot g(x,y)=F(x_a,y_b) +\varepsilon \cdot G(x_a,y_b)
     \end{displaymath}
     with $F(x_a,y_b),G(x_a,y_b)\in \langle x_a,y_b\rangle^m$.
     Again, by Taylor expansion and since $\varepsilon^2=0$ we have
     \begin{displaymath}
       f(x,y)=
       f(x_a,y_b) -
       \varepsilon\cdot \left(a\cdot {\textstyle\frac{\partial f}{\partial
           x}}(x_a,y_b)
         +b\cdot {\textstyle\frac{\partial f}{\partial y}}(x_a,y_b)\right)
     \end{displaymath}
     and
     \begin{displaymath}
       \varepsilon\cdot g(x,y)=\varepsilon\cdot g(x_a,y_b).
     \end{displaymath}
     Thus
     \begin{displaymath}
       F(x_a,y_b)=f(x_a,y_b)
     \end{displaymath}
     and
     \begin{displaymath}
       \langle x_a,y_b\rangle^m\ni G(x_a,y_b)=g(x_a,y_b)
       -a\cdot
       {\textstyle\frac{\partial f}{\partial
           x}}(x_a,y_b)
         -b\cdot {\textstyle\frac{\partial f}{\partial y}}(x_a,y_b).
     \end{displaymath}
   \end{proof}

   \begin{example}
     If we fix a curve $C\subset S$ and a point $p\in C$ such that
     $\mult_p(C)=m$, i.e.\ if using the notation of Example
     \ref{ex:incidencevariety}  we fix a point $C_p=(C,p)\in{\kh_m}$, then the
     diagram
     \begin{equation}\label{eq:def:4}
       \xymatrix{
         S\ar[r]^(0.4)\cong&\{C_p\}\times S\ar@{^{(}->}[rr]&&{\kh_m}\times S\ar[ddl]\\
         C\ar[r]^(0.25)\cong\ar[d]\ar@{^{(}->}[u]
         &\{(C_p,q)\;|\;q\in C\}\ar@{^{(}->}[r]\ar[d]\ar@{^{(}->}[u]
         &{\kf_m}\;\;\ar@{^{(}->}[ur]\ar[d]\\
         t_0\ar@{=}[r]&C_p\ar@{|->}[r]&{\kh_m}\ar@(ul,dl)[u]_\sigma
       }
     \end{equation}
     is an embedded equimultiple deformation of $C$ in $S$ along the
     section $\sigma$ through $p$. Moreover, any embedded equimultiple
     deformation of $C$ in $S$ with section through $p$ as a
     family is up to isomorphism induced via \eqref{eq:def:1} in a
     unique way  and thus
     factors obviously uniquely through \eqref{eq:def:4}. This
     means that every equimultiple deformation of $C$ in $S$ through
     $p$ is induced up to isomorphism in a unique way from
     \eqref{eq:def:4}.

     We now want to examine the tangent space to ${\kh_m}$ at a point
     $C_p=(C,p)$, which is just
     \begin{displaymath}
       T_{C_p}({\kh_m})=\Hom_{loc-K-Alg}\big(\ko_{{\kh_m},{C_p}},\C[\varepsilon]\big)
       =\Hom\big(T_\varepsilon,({\kh_m},{C_p})\big),
     \end{displaymath}
     where $({\kh_m},{C_p})$ denotes the germ of ${\kh_m}$ at ${C_p}$. However, a morphism
     \begin{displaymath}
       \psi:T_\varepsilon\longrightarrow ({\kh_m},{C_p})
     \end{displaymath}
     gives rise to a commutative fibre product diagram
     \begin{displaymath}
       \xymatrix{
         T_\varepsilon\times_{{\kh_m}}{\kf_m}\ar[r]\ar[d]^{\varphi'}&{\kf_m}\ar[d]\\
         T_\varepsilon\ar[r]^\psi\ar@(l,l)[u]^{\sigma'}&{\kh_m}\ar@(r,r)[u]_\sigma
       }
     \end{displaymath}
     sending the closed point of $T_\varepsilon$ to $C$.
     Thus $(\varphi',\sigma')\in \underline{\Def}^{sec,em}_{p\in C/S}(T_\varepsilon)$ is an
     \emph{embedded equimultiple deformation of $C$ in $S$ with
       section through $p$}. The
     universality of \eqref{eq:def:4} then implies that up to
     isomorphism each one
     is of this form for a unique $\varphi'$, and this construction
     is functorial. We thus have
     \begin{displaymath}
       T_{C_p}({\kh_m})\cong
       \underline{\Def}^{sec,em}_{p\in C/S}(T_\varepsilon),
     \end{displaymath}
     and hence \eqref{eq:def:3} gives the exact sequence
     \begin{displaymath}
       \xymatrix{
         0\ar[r] & \kk\ar[r] &
         T_{C_p}({\kh_m})\ar[r] &
         H^0\big(C,\kj_{Z/C}(C)\big)\ar[r] & 0.
       }
     \end{displaymath}
     In particular,
     \begin{displaymath}
       \dim_\C\big(T_{C_p}({\kh_m})\big)=
       \left\{
         \begin{array}{ll}
           \dim_\C H^0\big(C,\kj_{Z/C}(C)\big)-2, & \mbox{ if $C$ is
             unitangential},\\
           \dim_\C H^0\big(C,\kj_{Z/C}(C)\big)-1, & \mbox{ else}.
         \end{array}
       \right.
     \end{displaymath}
   \end{example}

   \begin{example} \label{ex:tangentspace}
     If we do the same constructions replacing in \eqref{eq:def:4} the
     family \eqref{eq:def:1a} by \eqref{eq:def:1als} we get for the
     tangent space to $\kl_m$ at $C_p=(C,p)$ the diagram of exact
     sequences
     \begin{displaymath}
       \xymatrix{
         0\ar[r] & \kk\ar[r] &
         T_{C_p}(\kh_m)\ar[r] &
         H^0\big(C,\kj_{Z/C}(C)\big)\ar[r] & 0\\
         0\ar[r] & \kk\ar[r]\ar@{=}[u] &
         T_{C_p}({\kl_m})\ar[r]\ar@{^{(}->}[u] &
         H^0\big(S,\kj_Z(C)\big)/H^0(\ko_S)\ar[r]\ar@{^{(}->}[u] & 0.
       }
     \end{displaymath}
     In order to see this consider the exact sequence
     \begin{displaymath}
       0\rightarrow \ko_S\rightarrow \ko_S(C)\rightarrow
       \ko_C(C)\rightarrow 0
     \end{displaymath}
     induced from the structure sequence of $C$. This sequence shows
     that the tangent space to $|L|$ at $C$ considered as a subspace of
     the tangent space $H^0\big(C,\ko_C(C)\big)$ of $H$ at $C$ is just
     $H^0\big(S,\ko_S(C)\big)/H^0(S,\ko_S)$ -- that is, a global
     section of $\ko_C(C)$ gives rise to an embedded deformation of
     $C$ in $S$ which is actually a deformation in the linear system
     $|L|$ if and only if it comes from a global section of
     $\ko_S(C)$, and the constant sections induce the trivial
     deformations. This construction carries over to the families
     \eqref{eq:def:1a} and \eqref{eq:def:1als}.
   \end{example}

     In particular we get the following proposition.

   \begin{proposition}\label{prop:expdim}
     Using the notation from above let $C$ be a curve in the linear
     system $|L|$ on $S$ and suppose that $p\in C$ such that
     $\mult_p(C)=m$.

     Then the tangent space of $\kl_m$ at
     $C_p=(C,p)$ satisfies
     \begin{displaymath}
       \dim_\C\big(T_{C_p}({\kl_m})\big)=
       \left\{
         \begin{array}{ll}
           \dim_\C H^0\big(S,\kj_Z(C)\big)-2, & \mbox{ if $C$ is
             unitangential},\\
           \dim_\C H^0\big(S,\kj_Z(C)\big)-1, & \mbox{ else}.
         \end{array}
       \right.
     \end{displaymath}
     Moreover, the expected dimension of $T_{C_p}(\kl_m)$ and thus of
     $\kl_m$ at $C_p$  is just
     \begin{displaymath}
       \expdim_{C_p}(\kl_m)=\expdim_\C\big(T_{C_p}({\kl_m})\big)=\dim|L|-\frac{(m+1)\cdot m}{2}+2.
     \end{displaymath}
   \end{proposition}
   For the last statement on the expected dimension just consider the
   exact sequence
   \begin{displaymath}
     0\rightarrow H^0\big(S,\kj_Z(C)\big)\rightarrow
     H^0\big(S,\ko_S(L)\big)\rightarrow H^0(S,\ko_Z)
   \end{displaymath}
   and note that the
   dimension of $H^0\big(S,\kj_Z(C)\big)$, and hence of $T_{C_p}(C)$,
   attains the  minimal possible value if
   the last map is surjective. The expected dimension of
   $H^0\big(S,\kj_Z(C)\big)$ hence is
   \begin{displaymath}
     \expdim_\C H^0\big(S,\kj_Z(C)\big)=\dim|L|+1-\deg(Z),
   \end{displaymath}
   and it suffices to calculate $\deg(Z)$. If $C$ is unitangential we
   may assume that $C$ locally in $p$ is given by $f=y^m+h.o.t.$, so that
   \begin{displaymath}
     \ko_{Z,p}=\C\{x,y\}/\langle y^{m-1}\rangle +\langle x,y\rangle^m,
   \end{displaymath}
   and hence $\deg(Z)=\frac{(m+1)\cdot m}{2}-1$. If $C$ is not
   unitangential, then we may assume that it locally in $p$ is given
   by an equation $f$ such that $f_m=\jet_m(f)=x^\mu\cdot y^\nu\cdot
   g$, where $x$ and $y$ do not divide
   $g$, but $\mu$ and $\nu$ are at least one. Suppose now that the
   partial derivatives of $f_m$  are not
   linearly independent, then we
   may assume $\frac{\partial f_m}{\partial x}\equiv\alpha\cdot
   \frac{\partial f_m}{\partial y}$ and thus
   \begin{displaymath}
     \mu y g\equiv \alpha \nu x g+
     \alpha x y\cdot \frac{\partial g}{\partial y}-
     x y\cdot \frac{\partial g}{\partial x},
   \end{displaymath}
   which would imply that $y$ divides $g$ in contradiction to our
   assumption. Thus the partial derivatives of $f_m$ are linearly
   independent, which shows that
   \begin{displaymath}
     \deg(Z)=\dim_\C\left(\C\{x,y\}/
       \langle {\textstyle\frac{\partial f}{\partial x},\frac{\partial
           f}{\partial y}}\rangle
       +\langle x,y\rangle^m\right)=\frac{(m+1)\cdot m}{2}-2.
   \end{displaymath}

   \begin{example}
     Let us consider the Example \ref{ex:incidencevariety-ls} in the case where $S=\pp^2$
     and $L=\ko_{\pp^2}(d)$. We will show that $\kl_m$ is then \emph{smooth of
     the expected dimension}. Note that $\pi(\kl_m)$ will only be smooth at $C$ if $C$ has an
     ordinary $m$-fold point, that is, if all tangents are different.

     \bigskip

     Given $C_p=(C,p)\in \kl_m$ we may pass to a suitable affine chart
     containing $p$ as origin and assume that
     $H^0\big(\pp^2,\ko_{\pp^2}(d)\big)$ is parametrised by
     polynomials
     \begin{displaymath}
       F_{\underline{a}}=f+\sum_{i+j=0}^d a_{i,j}\cdot x^iy^j,
     \end{displaymath}
     where $f$ is the equation of $C$ in this chart.
     The closure of $\pi(\kl_m)$ in $|L|$
     locally at $C$ is then given by several equations, say
     $F_1,\ldots,F_k\in \C[a_{i,j}|i+j=0,\ldots,d]$, in the
     coefficients $a_{i,j}$.  We get these equations by eliminating
     the variables $x$ and $y$ from the ideal defined by
     \begin{displaymath}
       \left\langle \frac{\partial^{i+j} F_{\underline{a}}}{\partial
           x^iy^j}
         \;\Big|\;i+j=0,\ldots,m-1\right\rangle.
     \end{displaymath}
     And $\kl_m$ is locally in $C_p$  described by the equations
     \begin{displaymath}
       F_1=0,\ldots,F_k=0,\;\;\frac{\partial^{i+j} F_{\underline{a}}}{\partial
           x^iy^j}=0,\;\;\;i+j=0,\ldots,m-1.
     \end{displaymath}
     However, the Jacoby matrix of these equations with respect to the
     variables $x,y,a_{i,j}$ contains a diagonal submatrix of size
     $\frac{m\cdot(m+1)}{2}$ with ones on the diagonal\tom{ -- we get
       these by differentiating the partial derivatives with respect
       to the $a_{i,j}$ --}, so that its
     rank is at least $\frac{m\cdot(m+1)}{2}$, which -- taking into
     account that $|L|=\pp\big(H^0\big(\pp^2,\ko_{\pp^2}(d)\big)\big)$
     -- implies that the
     tangent space to $\kl_m$ at $C_p$ has codimension at least
     $\frac{m\cdot(m+1)}{2}-1$ in the tangent space of $\kl$. By
     Proposition \ref{prop:expdim} we thus have
     \begin{align*}
       \dim_{C_p}(\kl_m) &\leq
       \dim_\C T_{C_p}(\kl_m)\leq \dim_\C
       T_{C_p}(\kl)-\frac{m\cdot(m+1)}{2}+1\\
       &=\dim(\kl)-\frac{m\cdot(m+1)}{2}+1\\
       &=\dim|L|-\frac{m\cdot(m+1)}{2}+2\\
       &=\expdim_{C_p}(\kl_m)\leq \dim_{C_p}(\kl_m),
     \end{align*}
     which shows that $\kl_m$ is smooth at $C_p$ of the expected dimension.
   \end{example}


\providecommand{\bysame}{\leavevmode\hbox to3em{\hrulefill}\thinspace}

\end{document}